\documentclass[12]{amsart}
\usepackage{amssymb,latexsym}
\usepackage[all]{xy}
\newcommand{\Ox}{{\mathcal O}}
\newcommand{\Px}{\mathbb P}
\newtheorem{thm}{Theorem}[section]
\newtheorem*{thm*}{Theorem}
\newtheorem{proposition}[thm]{Proposition}
\newtheorem{lemma}[thm]{Lemma}
\newtheorem{cor}[thm]{Corollary}
\theoremstyle{definition}
\newtheorem{defn}[thm]{Definition}
 
\newtheorem{notation}[thm]{Notation}  
\newtheorem{remark}[thm]{Remark}

\newcommand{\by}[1]{\stackrel{#1}{\longrightarrow}}

\newcommand{\boxtensor}{{\Box\kern-9.03pt\raise1.42pt\hbox{$\times$}}}
\newcommand{\F}{{\mathbb F}}
\newcommand{\Z}{{\mathbb Z}}  

\newcommand{\tensor}{\otimes}

\newcommand{\sF}{{\mathcal F}}

\newcommand{\sL}{{\mathcal L}}
\newcommand{\sM}{{\mathcal M}}

\newcommand{\sO}{{\mathcal O}}

\renewcommand{\tilde}{\widetilde}
 
\numberwithin{equation}{section}
\newcounter{elno}                

\newcounter{example}[section] 
\def\theexample{\thesection.\arabic{example}}

\begin{document}
\title{Semistability and Hilbert-Kunz multiplicities for 
curves}
{}

\author{V. Trivedi}
\address{School of Mathematics, Tata Institute of
Fundamental Research,
Homi Bhabha Road, Mumbai-400005, India}
\email{vija@math.tifr.res.in}

\date{}
\maketitle

\section{Introduction}
Let $(R,\bf{m})$ be a Noetherian local ring of
dimension
$d$ and 
of prime characteristic $p>0$,  and let $I$ be an ${\bf{m}}$-primary
ideal.
Then one defines the {\em Hilbert-Kunz function} of $R$ with respect
to $I$ as
\[HK_{R,I}(p^n)=\ell(R/I^{(p^n)}),\]
where~
\[I^{(p^n)}=\mbox{$n$-th Frobenius power of $I$}\]
\[=\mbox{ideal generated by $p^n$-th powers of elements of $I$}.\]
The  associated {\em Hilbert-Kunz multiplicity} is defined to be
\[HKM(R,I)=\lim_{n\rightarrow\infty} \frac{HK_{R,I}(p^n)}{p^{nd}}.\]

Similarly, for a non local ring $R$ (of prime characteristic $p$), and an
ideal $I\subseteq R$ for which $\ell(R/I)$ is finite, the Hilbert-Kunz
function and multiplicity make sense. Henceforth for such a pair $(R,
I)$, we denote the Hilbert-Kunz multiplicity of $R$ with respect to $I$ by
$HKM(R,I)$, or by $HKM(R)$ if $I$ happens to be an obvious maximal ideal. 

Given a pair $(X,\sL)$, where $X$ is a projective curve over an
algebraically closed field $k$ of positive characteristic $p$, and $\sL$
 is a base point free line bundle $\sL$ on $X$, define
$$HKM(X, \sL) = \mbox{HK multiplicity of the
section ring}~ B ~\mbox{with respect to the ideal}~ B_1B,$$
where $B = \oplus_{n\geq 0}H^0(X,\sL^{\tensor n})$ and  $B_1 =
H^0(X,\sL)$.
Note that when $\sL$ is very ample, giving an embedding $X\by{} 
\mathbf{P}^r_k$, then $HKM(X, \sL)$ equals 
the HK multiplicity of the ``homogeneous coordinate ring''$A = \oplus
A_n$, with respect to its maximal ideal $\oplus
A_{n > 0}$, where $A$ is  the image of  the natural map 
$\phi $, induced by $\sL$, 
$$\oplus_{n\geq 0} H^0(\mathbf{P}^r, \sO_{\mathbf{P}^r}(n))
\by{\phi}\oplus_{n\geq 0}H^0(X,\sL^{\tensor n}).$$ 
To discuss HK multiplicity of singular curves, we need to also
consider the  
HK multiplicity of $B$ with respect to the ideal generated by
$W\subseteq H^0(X,\sL)$, where $W$ is a base point free linear system,
which 
we denote by
$$HKM(X,\sL,W) = \mbox{HK multiplicity of}~~B~~\mbox{with respect to 
the ideal generated by}~~W . $$ 

\begin{notation}\label{n1} Now given $(X,\sL, W)$ as
above,
where $X$ is a nonsingular projective
curve
over $k$, consider the following  short exact sequence  
\begin{equation}\label{e1}
0 \by{} V_{\sL}(W) \by{} W\tensor\sO_X \by{}
\sL \by{} 0, \end{equation}
where $V_{\sL}(W)$ is a vector bundle of  rank $r =$ vector-space
dimension of $W - 1$ and is the 
kernel of the surjective map
$W\tensor \sO_X \by{} \sL$. If $W = H^0(X,\sL)$ then we denote
$V_{\sL}(W)$ by $V_{\sL}$.
\end{notation}
\vspace{5mm}

In Section~2, we prove (see Proposition~~\ref{2} and
Remark~\ref{r2}) 
that if $V_{\sL}$ is strongly semistable (\i.e., the pull back of 
$V_{\sL}$ under every iterated Frobenius map is semistable) then 
$$HKM(X, \sL) =
~\mbox{the HK multiplicity of the section ring with respect to its
graded maximal
ideal},$$
  (which may not be true in general without the strong semistability
condition). We also give a lower bound for
$HKM(X,\sL,W)$ in terms of $\deg~\sL$ and $\dim~W$, which is achieved when 
$V_{\sL}(W)$ is strongly semistable.
Later (see Theorem~\ref{16}) we prove 
 the
converse of this.

One consequence of Proposition~\ref{2} is that for given $(X,\sL)$, if
$HKM(X,\sL)$ does
not achieve the lower bound, then $V_{\sL}$ is not strongly
semistable.
 For a plane curve $X$ and $\sL = \sO_X(1)$, if $X$ is  nonsingular
or singular with certain
conditions on singularities 
then the referee
provided a 
proof (Proposition~\ref{p1}, corollaries~\ref{19} and \ref{ck}) that
$V_{\sL}$ is  semistable.

In Section~4, which has been rewritten as per the suggestions of
the referee,
we prove that, for an arbitrary base-point free ample line bundle $\sL$
on a nonsingular  curve $X$ of genus $g$ (hence for any irreducible
projective curve $C$), there is  an expression for
$HKM(X,\sL, W)$ (for $HKM(C, \sO_C(1))$) in terms of the ranks and
degrees of the vector bundles
occuring in a ``strongly stable Harder-Narasimhan filtration'' (in the
sense of recent work of A. Langer \cite{Langer}) of some Frobenius
pullback of $V_{\sL}(W)$ (see Theorem~\ref{sst}). Though this seems
difficult to 
use in actually computing the HK multiplicity, except when $
V_{\sL}(W)$ is
strongly
semistable, 
it does imply that it is a rational number, for instance. 
We also prove the converse to the Section~2 result mentioned above.

In Section~5, we discuss plane curves.  
In general, 
Theorem~\ref{18} gives a formula (and hence bounds) for the HK
multiplicity of an arbitrary 
plane curve $C$ of degree $d$ over a field of characteristic~$p
>0$.  In
particular
(Corollary~\ref{20}) if $X$
is a nonsingular plane curve of degree~$d$ then 
$$HKM(X,\sO_X(1)) =
\frac{3d}{4} + \frac{l^2}{4dp^{2s}}$$
where $0 \leq l\leq  d(d-3)$, and $l$ is an integer
congruent to $pd$ (mod 2), and $ s\geq 1$ (we allow $s = \infty $) is 
such that $F^{(s-1)*}V_{\sO_X(1)}$ is semistable and
$F^{s*}V_{\sO_X(1)}$
is not semistable (here $s=\infty$ means that $V_{\sO_X(1)}$ is strongly
semistable,).

The formulas (for singular and nonsingular plane curves) also imply that
for
$p>>0$ (for example when $p >
d(d-3)$), one can recover the numbers $s$ and $l$, where $l$ is the
measure of how much $F^{s*}V_{\sO_X(1)}$ is destablized, in the
sense that if $\sL_1 \subset F^{s*}V_{\sO_X(1)}$ is
the Harder-Narasimhan filtration then  slope~$\sL_1$ =
slope~$F^{s*}V_{\sO_X(1)} + l/2$. 
 So in this case, we have a
simple numerical
characterization of semistablity of the kernel bundle under the Frobenius
map via HK multiplicity.

 Using this, and Monsky's results ([M1], [M3]),  which are  explicit
computations
for certain nonsingular quartics), we prove the following 
(see Proposition~\ref{23}): for any integer $n\geq 1$, there exist 
explicit  rank $2$ vector bundles $V$ on nonsingular curves of genus $3$
over a field of characterstic $2$ or $3$, such that
$F^{n*}V$ is semistable, but $F^{(n+1)*}V$ is
not semistable. Morevover, when $ p=3$, the result also holds for $n=0$. 

 I would like to thank P. Monsky for his encouragement, as well as
interesting questions, and for sending me several of his as yet
unpublished papers. I would also like to thank R. Buchweitz for his
kind words. I would like to thank V. Srinivas for stimulating
discussions, and helpful suggestions.

Finally I would like to thank the referee for his detailed
comments and very useful suggestions and proofs.
  
Some of our results, particularly the formula for HK multiplicity in
Theorem~\ref{sst}, are also contained in an equivalent form in a recent
preprint of H. Brenner \cite{Brenner}. Our results here have been obtained
concurrently, and independently. The rationality of the HK multiplicity
of a smooth plane curve had  been also proved by Monsky (unpublished), by
different methods (private communications). 

\section{Semistability and HK multiplicity}

We first recall the notion of semistability. If $V$ is a vector bundle of
rank $r$ on a  
projective curve $X$, recall that ${\rm deg}~V:= {\rm deg}~(\wedge^r V)$,
and $\mbox{slope}~(V):= \mu(V) = {\rm deg}~V/{\rm rank}~V$. 

\begin{defn}\label{d1}Let $V$ be a vector bundle of rank $r$ on a
projective curve $X$. Then $V$ is {\it{semistable}} if for any 
subbundle
$V'\hookrightarrow V$, we have
  $$\mu(V') \leq \mu(V).$$ \end{defn}

\begin{defn}\label{d6}A vector bundle $V$ on $X$ is called {\em strongly
semistable}
if $F^{s*}V$ is semistable for the $s^{th}$ iterate of
the absolute Frobenius map,  $F^s: X\by{} X$, for all $s\geq 0$.\end{defn}

\begin{remark}\label{r1}If $W$ is
a line bundle then it
is semistable, and if $V$ is a semistable bundle
then so are $V^{\vee}$ and $V\tensor W$.
\end{remark}

From now onwards, $X$ is a nonsingular (projective) curve of genus $g\geq
2$ over an algebraically closed field $k$ of characteristic $p>0$ and
$\sL$ is a base point free line bundle on $X$, unless stated otherwise.
Recall the notation $h^i(X,\sF) := \dim_kH^i(X,\sF)$, for any coherent
sheaf $\sF$ on $X$, and $i=0,1$. 

\begin{lemma}\label{1}Let $X$ be a nonsingular projective curve of genus
$g$ and $V$ be a semistable bundle on $X$ of rank $r$ and degree $d$. Then
\begin{enumerate}
\item If ${\rm deg}~W <0$ then $h^0(X, W) = 0$,
\item If ${\rm deg}~W > r(2g-2)$ then $h^1(X, W) = 0$ and $h^0(X, W) =
{\rm deg}~W - r(g-1)$. 
\item If $0\leq {\rm deg}~W \leq r(2g-2)$ then $h^0(X,W) \leq rg$,
\end{enumerate}
\end{lemma}
\begin{proof} Statement (1) follows from the definition of semistable
vector bundle.
 
By Serre duality, we have  $h^1(X, W) = h^0(X,\omega_X\tensor
W^{\vee})$. Since $\omega_X\tensor W^{\vee}$
 is semistable, we get $h^0(X,\omega_X\tensor W^{\vee}) = 0$
if ${\rm deg}~W > r(2g-2)$, hence $h^1(X, W) =0$.
This, and the Riemann-Roch formula
$$h^0(X,W)- h^1(X,W) = \deg~W + r(1-g),$$
implies statement~(2). 

To prove statement~(3), we choose a line bundle $\sL$, given by
an effective divisor of degree~1, and an integer
$m\geq 0$ such that 
$\deg~(W\tensor\sL^m) \leq r(2g-2)$ and 
$\deg~(W\tensor\sL^{m+1}) > r(2g-2)$. Now 
$$h^0(X,W) \leq h^0(X, W\tensor\sL^{m+1}) = h^1(X,W\tensor \sL^{m+1}) +
\deg~(W\tensor\sL^{m+1})+r(1-g)$$
$$ = 
\deg~(W\tensor\sL^{m})+ r + r(1-g) \leq rg.$$
This proves statement~(3).\end{proof}

\begin{proposition}\label{2}Let $X$ be a nonsingular projective curve of 
genus $g$ and let $\sL$ be a base point free line bundle of degree $d$ on
$X$. If
$V_{\sL}$ (see (\ref{e1}))
is strongly semistable
then $$HKM(X,\sL) = HKM(B, {\bf m}) = 
 \frac{dh}{2(h-1)},$$
where $h = h^0(X,\sL)$, $B = \oplus_{n\geq 0} H^0(X,\sL^n)$ and 
 ${\bf{m}} = \oplus_{n> 0} H^0(X,\sL^n)$ is the graded maximal
ideal of $B$. 
\end{proposition}
\begin{proof}
 Let
$B_n =  H^0(X,\sL^n)$. Consider the Frobenius twisted multiplication map,
$$\mu_{k, n}:B_k^{(q)}\tensor B_{n-kq} \by{} B_n $$ 
given by 
$r\tensor r'\to r^qr'$, where $r\in B_k$ and $r'\in B_{n-kq}$ and
$B_k^{(q)} = B_k$ as an additive group  with $k$-action  on it given by 
$\lambda\cdot r = \lambda^q r$ for $\lambda \in k$ and $r\in B_k$.
Now 
$$\ell(B/{\bf{m}}^{(q)}) = \sum_n \ell(B_n/\sum_k \mbox{im}~\mu_{k,n}).$$
 Consider the short exact sequence 
$$0\by{} V_{\sL} \by{} H^0(X,\sL)\tensor \sO_X \by{} \sL \by{}
0,$$ 
 This gives 
$$0\by{} F^{s*}V_{\sL}\tensor \sL^{\tensor n} \by{}
H^0(X,\sL)^{(q)}\tensor \sL^{\tensor n} \by{} \sL^{\tensor n+q} \by{}
0,$$
where $q= p^s$ and $F:X\by{} X$ is the Frobenius map. 

Hence we have a
long exact sequence of cohomologies
$$H^0(X, F^{s*}V_{\sL}\tensor \sL^{\tensor n}) \by{}
H^0(X,\sL)^{(q)}\tensor H^0(X, \sL^{\tensor n})\by{} H^0(X,
\sL^{\tensor
n+q}) \by{}  H^1(X, F^{s*}V_{\sL}\tensor \sL^{\tensor n}),$$
where the second arrow is given by the map $\mu_{1,{n+q}}$.

Now ${\rm rank}~V_{\sL}
=
h-1$, and  
$$\begin{array}{lcl}
\deg~(F^{s*}V_{\sL}\tensor \sL^n) & = & \deg~(F^{s*}V_{\sL}) + (h-1)\deg
\sL^n  \\
 {} &  = & q \deg~V_{\sL} + (h-1) n (d)  \\
 {} & = & (-q + (h-1) n) d
\end{array}$$

\noindent~{\underline{Case}}~1\quad~~ Suppose
$n < q/(h-1) $. Then $\deg~(F^{s*}V_{\sL}\tensor \sL^n) < 0$. Hence by
Lemma~\ref{1}, the map $\mu_{1,n+q}$ is injective. 

Moreover
$n+q -kq < q/(h-1) + q - kq \leq 0$, if $k\geq 2$. In particular 
$\mbox{im}~ \mu_{k, n+q} = 0$ for $k\geq 2$. Hence in this range 
$\ell (B_{n+q}/\sum_k \mbox{im}~(\mu_{k, n+q})) = 
\ell (B_{n+q}/\mbox{im}~(\mu_{1, n+q})) = \ell(B_{n+q}) -
\ell(B_1)\cdot\ell(B_n)$.

\vspace{5mm}\noindent\quad~{\underline{Case}}~2\quad~~ Suppose $n >
q/(h-1)+ (2g-2)/d $. Then  $\deg~(F^{s*}V_{\sL}\tensor \sL^n) >
(h-1)(2g-2)$, hence by Lemma~\ref{1}, the map $\mu_{1,q}$ is surjective,
which implies
$\ell(B_{n+q}/\mbox{im}~(\mu_{1, n+q})) = 0$. Hence  
$\ell(B_{n+q}/\sum_k\mbox{im}~(\mu_{k, n+q})) = 0 $.

\vspace{5mm}\noindent\quad~{\underline{Case}}~3\quad~~Suppose $q/(h-1)
\leq n \leq q/(h-1) + (2g-2)/d$.
Then 
$$0 \leq \deg~(F^{s*}V_{\sL}\tensor\sL^n)
\leq (h-1)(2g-2),$$ and
therefore 

$$\sum_{n=\lfloor q/(h-1)\rfloor}^{
\lfloor q/(h-1)+(2g-2)/d\rfloor} h^0(X,
F^{s*}V_{\sL}\tensor
\sL^n) \leq (h-1)g\left(\frac{2g-2}{d} +1\right).$$

Therefore we have 
$$ HKM(X, \sL) = HKM(B,{\bf{m}})
 = \lim_{q\rightarrow\infty} \frac{1}{q^2} \sum_{n\geq
0}\ell(\frac{B_{n}}{{\rm im}{(\mu_{1,n})}})
 = \lim_{q\rightarrow\infty} \frac{1}{q^2} \sum_{n\geq
-q}\ell(\frac{B_{n}}{{\rm im}{(\mu_{1, n+q})}}) $$
$$ =  \lim_{q\rightarrow\infty} \frac{1}{q^2} \sum_{-q\leq n}
\left(h^0(X,\sL^{n+q}) -
h^0(X,\sL)h^0(X,\sL^{n}) + h^0(X,F^{s*}V_{\sL}\tensor\sL^n)\right)$$
$$ = \lim_{q\rightarrow\infty}\frac{1}{q^2}\sum_{-q\leq n
\leq q/(h-1)}h^0(X,\sL^{n+q}) -
h^0(X,\sL)h^0(X,\sL^{n})$$

$$ '' = \lim_{q\rightarrow\infty}\frac{1}{q^2}\sum_{0\leq n
\leq q/(h-1)+q}\chi(X,\sL^n) -h\sum_{0\leq n \leq
q/(h-1)}\chi(X,\sL^{n}) = (dh)/2(h-1)$$
This proves the proposition.\end{proof}

\begin{remark}\label{r2}In the above proof, replacing  
the complete
linear system by any base point free linear system $W$ of $\sL$, of
vector-space
dimension
$r+1$ (and replacing $h$ by $r+1$ everywhere), one sees that  if
$V_{\sL}(W)$ is
strongly
semistable then 
$HKM(X,\sL, W)= d(r+1)/2r$.
\end{remark}

  \section{Applications
and examples} In this section $X$ is a nonsingular curve and $\sL$ is a
base point free line bundle on $X$, and $V_{\sL}$ is the kernel vector
bundle
given
by the natural map
$$0\by{} V_{\sL}\by{} H^0(X,\sL)\tensor \sO_X \by{} \sL\by{} 0.$$

We use the following notation in this and in the forthcoming
sections.

\begin{notation}\label{n3}$C$ denotes an irreducible
curve of degree
$d >1$, over an
algebraically closed field
of characteristic $p$ and $\pi:X_C\by{} C$ is  the normalization of $C$,
where 
$g$ is  the genus of $X_C$ and  ${\sL_C} = \pi^*\sO_C(1)$ and  $W_C =
H^0(C,\sO_C(1))$. Note that  $W_C\subset H^0(X_C,{\sL_C})$ is a base point
free
linear
system. Hence this gives a natural short exact sequence 
of $\sO_{X_C}$-modules
\begin{equation}\label{e2}
0\by{} V_C \by{} W_C\tensor \sO_{X_C} \by{}
{\sL_C} \by{} 0,
\end{equation}
 where 
$V_C= V_{\sL_C}(W_C)$ following our earlier Notation~\ref{e1}.
\end{notation}

\begin{remark}\label{rn3}
Since $\pi$ is a finite birational map,
 by lemma~1.3 in [M0], theorem~2.7 in
[WY] or in [BCP], we have 
$$HKM(C,\sO_C(1)) = HKM(X_C,{\sL_C}, W_C).$$
 \end{remark}

Here we discuss
some examples $(X,\sL)$ for which the vector bundle $V_{\sL}$ is strongly
semistable.
 But before that we need to check the first necessary condition,
i.e., that the vector bundle $V_{\sL}$ is itself semistable.
The referee has provided the proofs of  Proposition~\ref{p1} and its
Corollaries~\ref{19} and ~\ref{ck}. Before coming to that 
we recall the following definition.

\vspace{5pt}

\begin{defn}\label{d0} The {\em gonality} of a nonsingular curve
$X$ is  the least
integer
$d$, for which there exists a line bundle of degree 
$d$ with a base point free
complete linear system of projective dimension 1 (in other words a line
bundle of
degree $d$ which induces a nonconstant map $X\to {\bf P}^1$).
\end{defn}

\begin{proposition}\label{p1}If 
 $X_C$ has gonality $\geq d/2$ then $V_{\sL}$ is
semistable.\end{proposition}
 \begin{proof} If $V_{\sL}$ is not
semistable,
then neither is 
$V_{\sL}^{\vee}$. Hence there exists a quotient line bundle $\sL_1$ of 
$V_{\sL}^{\vee}$
such that
$\mu(\sL_1) <
\mu(V_{\sL}^{\vee}) = d/2$. Since $V_{\sL}^{\vee}$ is globally generated,
the line
bundle $\sL_1$ is globally generated. Now $\sL_1$ cannot be the trivial
bundle; otherwise we will have $\sO_X \hookrightarrow V_{\sL}$ which
would imply
that 
$H^0(X,V_{\sL}) \neq 0$. So $h^0(X,\sL_1) \geq 2$. So it follows that $X$
has a line bundle, of degree $ < d/2$, with a
 linear system  of vector-space dimension $\geq 2$,
hence
a line
bundle of degree $< d/2$ with a base point free complete linear system 
of vector space dimension~2. In other words
the gonality of $X < d/2$, which  contradicts the hypothesis. This proves
the proposition.\end{proof}

\begin{cor}\label{19}If $X$ is a nonsingular plane curve, then
$V_{\sL}$, where $\sL= \sO_X(1)$,
is semistable.\end{cor}
\begin{proof}A classical result of M. Noether
(see [H] ,
theorem~2.1) implies that the gonality of $X$ is $d-1$, where $d$ is the
degree of $X$.  Now the proof follows from  Proposition~\ref{p1}.
\end{proof}

\begin{cor}\label{ck}Suppose $C$ is an irreducible projective plane curve
of
degree $d$ such
that the only singularities of $C$ are nodes and cusps, that $d\geq 4$ and
the number of singularities, $\delta$, satisfies $1\leq \delta \leq
d-2$. Then $V_C$ is semistable.\end{cor}
\begin{proof}Theorem~2.1 of [CK] implies (for $k=1$ in their
notation) that the gonality of $X_C$ is $\geq d-2$. Hence once again the
proof
follows
from Proposition~\ref{p1}.\end{proof}

In this context, we would also like to recall the following result
given in [T1],
which
was the  main ingredient in proving a conjecture of Monsky (see
Remark~\ref{r6} of this paper).

\begin{proposition}\label{p2}Let $C$ be an irreducible
projective plane curve of degree $d$ with a singularity of multiplicity
 $r \geq d/2$. Then: 
\begin{enumerate}
\item if $r = d/2$ then $V_C$ is strongly
semistable,
\item if $r > d/2$ then $V_C$ is not semistable and its
destabilizing line bundle is of degree $= r-d $.
\end{enumerate}
\end{proposition}

\section{HK multiplicities for base point free line bundles}
In this section, we consider $HKM(X,\sL,W)$ where $X$ is any non-singular
projective curve of genus $g$ over an algebraically closed field $k$ of
characteristic $p>0$, and $\sL$ is a line bundle on $X$ of degree $d$ with
base point
free linear system $W$.  
We derive an expression for the HK multiplicity in this case, involving 
terms which seem to be very difficult to compute, but which shows
that it is a rational number, with a denominator of a particular form.
As a consequence (see Remark~\ref{rn3}) the rationality of the HK
multiplicity of an
irreducible projective curve follows.

As mentioned in the introduction, this result  was
obtained independently by H. Brenner [B].  The
tools, both in Brenner's proof and ours, are Lemma~\ref{1},
Lemma~\ref{11}, and a recent result
of A. Langer [L] (Theorem~\ref{12}).  We shall also give a converse
to our Remark~\ref{r2}.

\begin{defn}\label{d2}Given a vector bundle  $E$ on $X$,  
a filtration by vector subbundles 
\[0=E_0\subset E_1\subset\cdots\subset
E_t\subset
E_{t+1}=E\]
is called a {\em Harder-Narasimhan filtration} (HN filtration) if 
\begin{enumerate}
\item[(i)]$E_1, E_2/E_1, \ldots, E_{t+1}/E_t$ are semistable vector
bundles, 
\item[(ii)] $\mu(E_1)>\mu(E_2/E_1)>\cdots>\mu(E_{t+1}/E_t)$.
\end{enumerate}
\end{defn}

\begin{remark}\label{r3}Note that such a filtration 
exists and
is
unique (see [HN],
lemma~1.3.7).  Moreover, if $t\geq 1$, then  
$$\mu(E_i) > \mu(E_i/E_{i-1}),~~\mbox{for all}~~ 2\leq i\leq t+1.$$ 
The case when $E$ is semistable corresponds to $t =0$.
\end{remark}

\vspace{5pt}

\begin{notation}\label{n4}~~ If $ 0\subset
E_1\subset\cdots\subset
E_t\subset
E_{t+1}= E$ is the HN filtration of $E$ then we write
$$\mu_{{\rm max}}(E) = \mu(E_1)~~\mbox{and}~~\mu_{{\rm min}}(E) =
\mu(E/E_t).$$
\end{notation}

\begin{defn}\label{d3} A filtration of subbundles
\[0=E_0\subset E_1\subset\cdots\subset E_t\subset E_{t+1}=E\] 
of $E$ is a {\em strongly stable HN
filtration} if it is a HN filtration and 
$E_1, E_2/E_1, \ldots, E_{t+1}/E_t$ are strongly semistable vector
bundles. \end{defn}

Note that whenever $E$ has a strongly stable HN
filtration
then the
HN-filtration of $F^{k*}(E)$ is 
$$ 0 \subset F^{k*}(E_1)\subset F^{k*}(E_2)\subset \cdots \subset
F^{k*}(E_t)\subset
F^{k*}(E_{t+1}) =F^{k*}(E).$$

Now recall  the crucial result of Langer~[L], which we state for the
special case of curves.

\vspace{5mm}

\begin{thm}~\label{12}{(A. Langer)}~\quad If $V$ is a vector
bundle on a
nonsingular projective curve defined over an algebraically closed
field of characteristic $p > 0$, then there exist $s>0$ such that 
$F^{s*}(V)$ has a strongly stable HN filtration.\end{thm}

\begin{defn}\label{d4}For a vector bundle $V$ on $X$, and an ample line
bundle
$\sL$ on $X$, we
define 
$$\sigma_s(V)=\sum_{n\le 0}h^0(F^{s*}(V)\otimes \sL^n)+
\sum_{n> 0}h^1(F^{s*}(V)\otimes \sL^n).$$
\end{defn}

\begin{lemma}\label{13}
If $V$ is a strongly semistable vector bundle of rank $r$ and degree $a$,
and $\deg~\sL = d$,
then 
$$\sigma_s(V)=\frac{a^2}{2rd}p^{2s}+O(p^s).$$
\end{lemma}
\begin{proof} Suppose for example that $a\ge 0$.  We are given that 
$F^{s*}(V)\otimes \sL^n$ is semistable of degree
$p^sa+rdn$. We choose $s>0$ such that $(2g-2)/d < p^sa/rd$. Then
$$\sigma_s(V) =
\displaystyle{\sum_{n<\frac{-p^sa}{rd}}}h^0(X,F^{s*}(V)\otimes \sL^n) +
\displaystyle{\sum_{\frac{-p^sa}{rd}\leq n\leq
\frac{2g-2}{d}-\frac{p^sa}{rd}}}
h^0(X,F^{s*}(V)\otimes \sL^n) $$

$$ +\displaystyle{\sum_{\frac{2g-2}{d}-\frac{p^sa}{rd} < n \leq
0}}h^0(X,F^{s*}(V)\otimes
 \sL^n)
+ \sum_{n >0} h^1(X,F^{s*}(V)\otimes \sL^n).$$
Now applying Lemma~\ref{1} to this equation we get 
$$\sigma_s(V) = C_0 +\sum_{\frac{2g-2}{d}-\frac{p^sa}{rd} < n 
\leq 0}h^0(X,F^{s*}(V)\otimes \sL^n) 
= C_0 +\sum_{\frac{2g-2}{d}-\frac{p^sa}{rd} < n 
\leq 0}\chi(X,F^{s*}(V)\otimes \sL^n),$$
where $0\leq C_0\leq rg((2g-2)/d +1)$. This gives 
$\sigma_s(V)=\frac{a^2}{2rd}p^{2s}+O(p^s)$. 
The argument for $a< 0$ is similar.\end{proof}

\begin{notation}\label{n5}To generalize Lemma~\ref{13} to
an arbitrary vector bundle $V$ on $X$,
we shall attach a rational number $\alpha(V)$ to $V$, as follows.
We choose  $m\geq 0$ such that the vector bundle $F^{m*}V$
has a 
strongly stable HN filtration (this is possible by Theorem~\ref{12}),
$$0\subset E_1\subset E_2\subset \cdots \subset E_t \subset E_{t+1} =
F^{m*}V,$$
Recall that, for any $n\geq 0$,
$$0\subset F^{n*}E_1\subset F^{n*}E_2\subset \cdots \subset F^{n*}E_t
\subset
 F^{n*}E_{t+1} = F^{(m+n)*}V,$$
is the strongly stable HN filtration of $F^{(m+n)*}V$.
We set 
$$a_i=p^{-m}\mbox{deg}(E_i/E_{i-1}),
~~r_i=\mbox{rank}(E_i/E_{i-1})$$
\begin{equation}\label{e3}
\alpha(V)=\sum_i(a_i^2/r_i).
\end{equation}
\end{notation}

\begin{remark}\label{r4}Note
that these numbers are  independent of the choice of $m$, and 
that 
$$\sum a_i=a,~~\mbox{and }~~ \sum r_i=r.$$
\end{remark}

\begin{lemma}\label{11}
Let
$0\to U\to V\to W\to 0$
be an exact sequence
of vector bundles on
$X$. Suppose that $U$ and $V$ admit strongly stable HN filtrations, and
that 
$$\mu_{\rm min}(U)-\mu_{\rm max}(W) > {\rm max}(0,2g-2).$$ Then
$\sigma_s(V)=\sigma_s(U)+\sigma_s(W)$ for all $s$.
\end{lemma}
\begin{proof}
It suffices to show that 
\[h^0(X,F^{s*}(V)\tensor\sL^n)=h^0(X,F^{s*}(U)\tensor\sL^n)+
h^0(X,F^{s*}(W)\tensor \sL^n)\]
for all $s$ and $n$.
Consider the canonical long exact sequence
$$0\by{} H^0(F^{s*}(U)\tensor\sL^n) \by{} H^0(F^{s*}(V)\tensor\sL^n) \by{}
H^0(F^{s*}(W)\tensor\sL^n) \by{} H^1(F^{s*}(U)\tensor\sL^n) \by{} .$$  
Now 
$$\mu_{\rm min}(F^{s*}(U)\tensor \sL^n) -\mu_{\rm max}(F^{s*}(W)\tensor
\sL^n) = p^s(\mu_{\rm min}(U)-\mu_{\rm max}(W) > 2g-2.$$
Therefore, either $\mu_{\rm max}(F^{s*}(W)\tensor \sL^n) < 0$, in which
case 
$h^0(F^{s*}(W)\tensor \sL^n) = 0$, or 
$$\mu_{\rm min}(F^{s*}(U)\tensor \sL^n) > 2g-2,$$
 in which case,  we have  
$h^1(F^{s*}(U)\tensor \sL^n) = 0$, by Serre duality. 
Hence the lemma follows, by the above long exact sequence.\end{proof}

\begin{cor}\label{14}For any vector-bundle $V$ on $X$, 
$$ \sigma_s(V)=\frac{\alpha(V)}
{2d}p^{2s}+O(p^s).$$\end{cor}
\begin{proof}Taking large enough Frobenius pull backs, {\em i.e.} for
$ m >> 0$, we can make sure
that
$$0\subset E_1\subset E_2\subset \cdots \subset E_t \subset E_{t+1} =
F^{m*}V $$
is the strongly stable HN filtration of $F^{m*}V$ and
$$\mu(E_i/E_{i-1}) -  \mu(E_{i+1}/E_i) > r(2g-2),$$
hence, by Remark~\ref{r3},
$$\mu(E_i) -  \mu(E_{i+1}/E_i) > r(2g-2).$$
Moreover, $E_{i+1}/E_i$ is strongly semistable  and $0\subset E_1\subset
\cdots \subset E_i$ is the strongly stable HN filtration of $E_i$. 
Hence applying Lemma~\ref{11}, for $s-m >0$ we get 
$$\sigma_{s-m}(E_{i+1}) = \sigma_{s-m}(E_i) + \sigma_{s-m}(E_{i+1}/E_i).$$
Now, for $s-m>>0$, by induction
$$\sigma_{s}(V) = \sigma_{s-m}(E_{t+1}) = \sigma_{s-m}(E_{1}) + 
\sigma_{s-m}(E_{2}/E_1) +\cdots
\sigma_{s-m}(E_{t+1}/E_t).$$  
Now the corollary follows from Lemma~\ref{13}.
 \end{proof}

\begin{thm}\label{sst}Let $X\subset \Px^r$ be a nonsingular projective
curve over $k$
and let $\sL$ be a 
line bundle on $X$ of degree $d$, with a base point free linear system
$W$.
Then 
$$HKM(X,\sL,W) = (1/2d)(d^2+\alpha(V_{\sL}(W))).$$
In particular $HKM(X, \sL,W)$ is a rational number.\end{thm}
\begin{proof}Let $B$ be the section ring
$\oplus_{n\geq 0} H^0(X, \sL^n)$, and $I$ be the ideal
of $B$ generated by $W\cdot B$.
We only need show that the HK multiplicity
of $B$ with respect to $I$ is 
$(1/2d)(d^2+\alpha(V_{\sL}(W)))$.  Making use of the
various exact sequences
$$0\to F^{s*}(V_{\sL}(W))\otimes \sL^n\to
\sL^n\oplus\cdots \oplus \sL^n\to 
\sL^{n+p^s}\to 0,$$
one finds that 
$$\dim~\frac{B}{I^{[p^s]}B} = \sum_n(h^0(X,F^{s*}(V_{\sL}(W))\otimes
\sL^n)-
(r+1)h^0(X,\sL^n)+h^0(X,\sL^{n+p^s})).$$

Now each term in this sum is unchanged when
$h^0$ is replaced by $h^1$.  So the sum is
$$\sigma_s(V_{\sL}(W))-(r+1)\sigma_s(\Ox_X)+\sigma_s(\sL).$$
Since $\alpha(\Ox_X)=0$ and $\alpha(\sL)=d^2$,
by Corollary~\ref{14}, we have
$$\mbox{dim}(B/I^{[p^s]}B) = \frac{1}{2d}(\alpha(V_{\sL}(W))+d^2)p^{2s}+O(p^s).$$
This proves the theorem.\end{proof}

\begin{remark}\label{r5}We have
$$\frac{b^2}{s} + \frac{c^2}{t} - \frac{(b+c)^2}{s+t}
 = \frac{(cs-bt)^2}{st(s+t)}.$$
So if $s,t>0$, 
$$\frac{b^2}{s} + \frac{c^2} {t} 
\ge \frac{(b+c)^2}{s+t},$$
with equality if and only if $b/s=c/t$. It follows
that $\alpha(V_{\sL}(W))\ge d^2/r$ with equality if
and only if $V_{\sL}(W)$ is strongly semistable.
Together with Theorem~\ref{sst}, this gives: 
\end{remark}

\begin{thm}\label{16}For a nonsingular projective curve $X$ 
with a line bundle $\sL$ of degree $d$ and a base point free linear system
$W$, of
$\sL$, of 
dimension $r$, 
$$HKM(X,\sL,W) \geq d(r+1)/2r ,$$
and 
$$HKM(X,\sL,W) =  d(r+1)/2r$$
if and only if $V_{\sL}(W)$ is
strongly semistable. 
 \end{thm}

Now, Remark~\ref{rn3} implies the following

\begin{cor}\label{cc}If $C\subseteq {\mathbf P}^r$ is an irreducible
projective curve of degree
$d$
then 
$$HKM(C,\sO_C(1)) = (1/2d)(d^2+\alpha(V_C)),$$
which is a rational number.
Furthermore 
$$HKM(C,\sO_C(1)) \geq d(r+1)/2r,$$
with equality  if and only if
 $V_C$ is
strongly semistable\end{cor}

\begin{cor}\label{c1}If $X$ is a nonsingular projective curve of genus
$g\geq 2$ and 
$\omega_X$ is the canonical sheaf of $X$ then
$$HKM(X,\omega_X)\geq g, $$
with equality if and only if 
$V_{\omega_X}$ is stongly semistable.\end{cor}

\section{HK multiplicity for plane curves}
In this section we use the Notation~\ref{n3},
 where
$C$ is
an  irreducible
plane curve of degree $d >1$, over an
algebraically closed field
of characteristic $p$.
 Hence we have a natural short exact sequence 
of $\sO_{X_C}$-modules
$$0\by{} V_C \by{} W\tensor \sO_{X_C} \by{} {\sL_C} \by{} 0,$$
where $V_C = V_{\sL}(W)$ is a rank two vector bundle.

\vspace{5pt}

\begin{remark}~\label{r55}For a rank two vector bundle
$V$,
either the
bundle
is
strongly semistable or some iterated Frobenius pull back has  HN
filtration given by a line
bundle ${\sL} \subset F^{s*}V$ such that $F^{s*}V/{\sL}$ is also a line
bundle. In other
words the HN filtration of $F^{s*}V$ is a strongly stable HN
filtration. Hence the result
of Langer is obvious.
\end{remark}
\vspace{5mm}

The following lemma is proved in [SB], Corollary~$2^p$ (see also [L]). 
We sketch another proof.

\begin{lemma}\label{17}Let $X$ be a nonsingular curve of genus $g$ over
an
algebraically closed field $k$ of characteristic $p>0$. Let $V$ be a
 vector bundle of rank 2
over $X$. Suppose there exists an exact sequence  
$$ 0\to \sL_1 \to F^*V \to \sM_1 \to 0, $$ 
such that $\sL_1$, $\sM_1$ are line bundles, and
$$\deg~\sL_1-\deg~\sM_1 > {\rm max}~(2g-2, 0).$$ Then $V$ is not
semistable.\end{lemma}
\begin{proof}If $g=0$ and $V$ is semistable then $F^*(V)$ is
semistable. This contradicts the  hypothesis that 
$\deg~\sL_1 - \deg~\sM_1 > 0$. So we may assume that $g>0$. Hence
$\deg~\sL_1 - \deg~\sM_1 > 2g-2 $.
Then there is
a
canonical connection 
$\nabla:F^*(V) \by{} F^*(V)\tensor\omega_X$ given locally by 
$$\nabla(F^*(e_1)) = \nabla(F^*(e_2)) = 0,$$
where $\{e_1,e_2\}$ is any local basis for $V$. Let $f = p\circ
\nabla\mid_{\sL_1}$, where $p:F^*(V)\tensor \omega_X \by{}
\sM_1\tensor\omega_X$ is the obvious map. Let $a$ and $s$ be local
sections
of $\sO_X$ and $\sL_1$ respectively. Then 
$$f(as) = p(s\tensor da + a\nabla s) = p(a \nabla s) = a f(s).$$
Hence $f:\sL_1\by{} \sM_1\tensor \omega_X$ is an $\sO_X$-linear map. 

If $f\neq 0$ then $\deg~\sL_1
 \leq \deg~\sM_1 + (2g-2)$ which would
contradict
the hypothesis. Hence $f =0$.  Now, note that
locally, $\sL_1$ is a free $\sO_X$-module of rank 1 in $F^*V$, generated
by a section  of the form  $s = aF^* e_1 + F^* e_2$, or of the form $s =
F^* e_1 + bF^* e_2$. 
Without loss of generality one can assume $s = aF^* e_1 + F^* e_2$. Then 
$f(s)  = 0$ implies $F^* e_1 \tensor da \in \sL_1\tensor \omega_X$. Hence
we
can find a local section $w$ of $\omega_X$ such that 
$ F^*e_1\tensor da =  (aF^* e_1 + F^* e_2) \tensor w $, which implies
$w=0$ and $da= 0$. Hence $a = \tilde{a}^p$ for some local section
$\tilde{a}$ of $\sO_X$. This implies $aF^*e_1 + F^* e_2 = F^*(\tilde{a}e_1
+ e_2)$.  Hence $\sL_1 = F^*{\sL_1'}$ for some line sub-bundle $\sL_1'$ of
$V$.
Since  $\deg~F^*(\sL_1') > 1/2\deg~F^*(V)$ we have $\deg~\sL_1' > \mu(W)$,
which implies that $V$ is not semistable.
\end{proof}

\begin{thm}\label{18}Let $C$ be an irreducible plane curve of degree $d
>1$. Let $X_C \by{\pi} C$ be the normalization of $C$. Let $V_C$ be the
rank
two vector bundle given by the natural map
$$0\by{} V_C \by{} H^0(C,\sO_C(1))\tensor\sO_X \by{} {\sL_C} \by{}
0.$$
Then one of the following holds:
\begin{enumerate}
\item  $V_C$ is strongly semistable. In this case 
$HKM(C) = 3d/4$.
\item $V_C$ is not  semistable. Then 
$$HKM(C) = \frac{3d}{4} + \frac{l^2}{4d},$$
where $0 < l <d$ and $l$ is an integer congruent to $d$ (mod 2). 
\item $V_C$ is semistable but not strongly semistable. Let $s \geq 1$
be
the number such that $F^{(s-1)*}V_C$ is semistable and $F^{s*}V_C$ is not
semistable. Then 
$$HKM(C) = \frac{3d}{4} + \frac{l^2}{4dp^{2s}},$$
where $l$ is an integer  congruent to $pd$ (mod 2) with 
$0 < l \leq 2g-2 $, so that in particular $0<l \leq d(d-3)$.
\end{enumerate}
\end{thm}
\begin{proof}{\underline{Case}}~(1) follows from Remark~\ref{r2} with
$r=2$.

\noindent{\underline{Case}}~(2) Given that $V_C$ is not
semistable, we have 
\[0\to\sL_1\to V_C \to\sM_1\to 0\] 
where 
$$\mu(\sL_1) = \deg~\sL_1 = -\frac{d}{2} +\frac{l}{2}
~~~~~~\mbox{and}~~~~~~
\mu(\sM_1) = \deg~\sM_1 = -\frac{d}{2} - \frac{l}{2},$$
for some $l>0$ and $l$ is an integer congruent to $d (mod~2)$.
Since this is the strongly stable HN filtration (see Remark~\ref{r55}), by
Theorem~\ref{sst}
$$HKM(C) = \frac{3d}{4} + \frac{l^2}{4d} .$$
Since an irreducible plane curve of degree $d>1$ has HK multiplicity $<d$,
we have $l <d$. This proves the statement~(2).

\noindent{\underline{ Case}}~(3). 
If $\sL_1$ is the destabilizing bundle of $F^{s*}V_C$ then there  exists a
short
exact sequence 
$$0\by{} \sL_1 \by{} F^{s*}V_C \by{} \sM_1 \by{} 0, $$ \
such that for some positive integer $l$
$$\deg~\sM_1 = -\frac{d}{2}p^{s} - l/2,~~~\mbox{and}~~ 
\deg~\sL_1 = -\frac{d}{2}p^{s}+ l/2.$$
Since $F^{(s-1)*}V_C$ is semistable, by Lemma~\ref{17}, we have
$$\deg~\sL_1 - \deg~\sM_1 = l 
\leq 2g-2 .$$
 Since $0\subset \sL_1 \subset F^{s*}V_C$ is the strongly stable HN
filtration,  Theorem~\ref{sst} and a calculation like that made in
case~(2) gives the desired value of $HKM(C)$. 
This proves the theorem.\end{proof}

If $X$ is a nonsingular plane curve, then by 
Corollary~\ref{19}, the bundle $V_{\sO_X(1)}$ is semistable,
and so Theorem~\ref{18} gives the following corollary.

\begin{cor}\label{20}Let $X$ be 
a nonsingular plane
curve of degree $d$ over an algebraically closed field of characteristic 
$p>0$, and 
$\sO_X(1)$ the corresponding very ample line bundle. 
Then $$
HKM(X,\sO_X(1))   =  \frac{3d}{4}+\frac{l^2}{4dp^{2s}},
$$
where $s\geq 1$ is a number such that $F^{(s-1)*}V_{\sO_X(1)}$ is
semistable
and 
$F^{s*}V_{\sO_X(1)}$ is not semistable (if $F^{t*}V_{\sO_X(1)}$ is 
semistable for all $t\geq 0$, we take
$s=\infty$) and $l$ is an integer congruent to 
$pd$ (mod 2) with 
$0\leq l \leq d(d-3)$.
\end{cor}

\begin{remark}\label{201}If all the singularities of an irreducible
projective 
plane curve of degree $d>1$ are  nodes and cusps,
and the number of singularities is  $\leq
d-2$, then, by Corollary~\ref{ck}, it follows that  Case~(2) of
Theorem~\ref{18} can not occur.\end{remark}

\begin{remark}\label{r6}Suppose $C$ is an irreducible
projective
plane curve with a
singularity of multiplicity
$ = r \geq d/2$. Monsky conjectured 
$$HKM(C) = \frac{3d}{4} + \frac{(2r-d)^2}{4d} .$$
We proved this in [T1]; note that it is an immediate consequence of
cases~(1) and (2) of Theorem~\ref{18}, combined with Proposition~\ref{p2}.
\end{remark}
\vspace{5pt}

\begin{remark}Let $C$ be an  irreducible plane quartic. If $C$ is
singular, the last remark shows that $HKM(C)$ is $3$ if $C$ has a point of
multiplicity $2$, and is $13/4$ if $C$ has a triple point. 

If  $C$ is nonsingular, then we are either in case~(1) of
Proposition~\ref{18}, or in case~(3) of the same proposition with
$l=2$
or $4$.  So $HKM(C)$ is either $3$, $3+(1/p^s)$ or $3+(1/4p^{2s})$, for
some $s\geq 1$. This result had been conjectured by Monsky.

In particular, when $C$ is nonsingular, we have  $HKM(C) \leq
3+(1/p^2)$. The
referee informs us that when $p =2$, we have 
 $HKM(C)\leq 3+ (1/16)$. 
\end{remark}

\vspace{5mm}

We recall some results  of Monsky~[M1], [M3] (see also [M2]), about 
nonsingular quartics of a certain type..

\begin{thm}\label{21}(Monsky) Let $R_{\alpha} = k[x,y,z]/(g_{\alpha})$,
where ${\mbox{char}~k} = 2$ and 
$$g_{\alpha} = \alpha x^2y^2 + z^4 +
xyz^2+(x^3+y^3)z,$$
with $\alpha \in k \setminus \{0\}$. Then 

$$HKM(R_{\alpha}) = 3 + 4^{-m(\alpha)},$$
where, for $\lambda \in k$ such that 
$\alpha = \lambda^2+\lambda$,  we define $m(\alpha)$ as follows:
$$\begin{array}{lcl}
m(\alpha) & = & \mbox{deg of $\lambda$ over $\Z/2\Z$ if
$\alpha$ is 
algebraic over}~ \Z/2\Z\\
  & = & \infty~ \mbox{if $\alpha$ is transcendental over}~ 
\Z/2\Z
\end{array}$$
\end{thm}

\begin{thm}\label{22}(Monsky) Let $R_{\lambda} =
k[x,y,z]/(f_{\lambda})$,
where ${\mbox{char}~k} = 3$
and
$$f_{\lambda} = z^4 - xy(x+y)(x+\lambda y),$$
with $\lambda \in k \setminus \{0, 1\}$. Then 

$$HKM(R_{\lambda}) = 3 + \frac{1}{p^{2d(\lambda)}},$$
where $d = d(\lambda)$ is the degree of $\lambda$ over 
$\Z/3\Z$ (and  $d = \infty$ if $\lambda$ is transcendental over 
$\Z/3\Z$).
\end{thm}

Note that $X_{\alpha} = \mbox{Proj}~R_{\alpha} \by{\pi} \mathbf{P}^2$ is a
nonsingular plane quartic of genus~3.  We also note that, 
 given any integer $n\geq 2$ there exists an $\alpha \in
\bar{\F_2}$ such that $m(\alpha) = n$. Similarly given any $n\geq 1$
there exists $\lambda \in \bar{\F_3}$ such that $d(\lambda) =
n$. 

Applying
Corollary~\ref{20} to  Example~\ref{21}, we see that
$F^{(n-1)*}V_{\alpha}$ is semistable and $F^{{n+1}*}V_{\alpha}$ is not. 
(The referee has shown that 
$F^{n*}V_{\alpha}$ is  semistable). Hence 
we get the following.

\begin{proposition}\label{23} \begin{enumerate}
\item[(i)] Given any integer $n\geq
2$, there exists a
non-singular quartic curve $X_{\alpha} \subseteq
\mathbf{P}^2_{\bar{\F_2}}$, given by the equation
\[\alpha x^2y^2 + z^4 + xyz^2+(x^3+y^3)z =0\]
where $m(\alpha)=n$, such that the vector bundle 
\[V_{\alpha} = \Omega^1_{\mathbf{P}^2}\mid_{X_{\alpha}}\] 
is a semistable vector bundle on $X_{\alpha}$ of rank 2 and degree -4,
and the iterated Frobenius pullback $F^{n*}V_{\alpha}$ is
not semistable,
while $F^{(n-1)*}V_{\alpha}$ is semistable.
\item[(ii)]
Given any integer $n\geq 1$, there exists a
non-singular quartic curve $X_{\lambda} \subseteq
\mathbf{P}^2_{\bar{\F_3}}$, given by the equation
\[z^4 - xy(x+y)(x+\lambda y) \]
where $d(\lambda)=n$, such that the vector bundle 
\[V_{\lambda} = \Omega^1_{\mathbf{P}^2}\mid_{X_{\lambda}}\] 
is a semistable vector bundle on $X_{\alpha}$ of rank 2 and degree -4,
and the iterated Frobenius pullback $F^{n*}V_{\lambda}$ is
not semistable,
while $F^{(n-1)*}V_{\lambda}$ is semistable.
\end{enumerate}
\end{proposition}

\begin{remark}\label{r11}Let $R_{\lambda}$ be as in Theorem~\ref{22}, but
with $p>3$. Monsky [M3] has given a practical algorithm involving the
iteration of a rational function, for calculating
$HKM(R_{\lambda})$. Together with our results, this lets one calculate the
smallest power of $F^*$ that destabilizes $V_{\lambda}$.
\end{remark}

\end{document}